\documentclass[11pt]{article}
\usepackage{amssymb}
\usepackage{multirow}
\usepackage{graphicx}
\usepackage{amsmath,amssymb,verbatim}
\usepackage{amsthm}
\usepackage{epsfig}
\usepackage[latin1]{inputenc}

\numberwithin{equation}{section}

\usepackage[main=english]{babel}
\newtheorem{definition}{Definition}[section]
\newtheorem{example}{Example}[section]
\newtheorem{theorem}{Theorem}[section]
\newtheorem{remark}{Remark}[section]
\newtheorem{lemma}{Lemma}[section]

\newcommand{\D}{\mathbb{D}}

\newcommand{\I}{\mathbb{I}}
\newcommand{\R}{\mathbb{R}}
\newcommand{\C}{\mathbb{C}}
\newcommand{\Com}{\mathbb{C}}

\begin{document}



\centerline{{\Large{\bf Symmetrical Sonin kernels 
}}}
\vspace{0.2cm}

\centerline{{\Large{\bf in terms of the hypergeometric functions
        }}}

\vspace{0.3cm}


\vspace{0.3cm}
\centerline{{\bf Yuri Luchko }}
\vspace{0.3cm}

\centerline{{Department of Mathematics, Physics, and Chemistry}}

\centerline{{Berlin University of Applied Sciences and Technology}}

\centerline{{Luxemburger Str. 10, 13353 Berlin,\ Germany}}

\centerline{e-mail: luchko@bht-berlin.de}

\begin{abstract}
In this paper, we introduce a new class of the kernels of the integral transforms of the Laplace convolution type that we call symmetrical Sonin kernels. For a symmetrical Sonin kernel given in terms of some elementary or special functions, its associated kernel has the same form with possibly different parameter values. Several known and new kernels of this type are derived by means of the Sonin method in the time domain and using the Laplace integral transform in the frequency domain. The new symmetrical Sonin kernels are provided in terms of the Wright function and some extensions of the Horn confluent hypergeometric functions in two variables. 
\end{abstract}


\vspace{0.2cm}

\noindent
{\sl MSC 2010}: 26A33 (primary);  33C60; 33C70; 44A05; 44A10

\noindent
{\sl Key Words}: Sonin kernels, symmetrical Sonin kernels,  Laplace integral transform, hypergeometric functions, Horn functions, general fractional derivative, general fractional integral



\section{Introduction}
\label{s1}

In the theory of the integral transforms of the Mellin convolution type, the functions $p$ and $q$ are said to form a pair of the Fourier kernels if the relations
\begin{equation}
\label{Four1}
F(x) = \int_0^{+\infty} p(x\xi)\, f(\xi)\, d\xi,\ x>0,
\end{equation}
\begin{equation}
\label{Four2}
f(\xi) = \int_0^{+\infty} q(\xi x)\, F(x)\, dx,\ \xi>0,
\end{equation}
are simultaneously valid on some functional spaces that the functions $f$ and $F$ belong to. In the case $p(x)\, = \, q(x),\ x\in \R_{+}$, the kernels $p$ and $q$ are called the symmetrical Fourier kernels. 

Let us suppose that the Mellin integral transforms 
\begin{equation}
\label{Mel}
P(s) = \int_0^{+\infty} p(x)\, x^{s-1}\, dx,\ \ Q(s) = \int_0^{+\infty} q(x)\, x^{s-1}\, dx
\end{equation}
of the Fourier kernels  $p$ and $q$  do exist on the domain $D_{(p,q)} = \{s \in \Com : c_{(p,q)}< \Re(s) < C_{(p,q)}$\}. The Mellin convolution theorem applied to 
the equations \eqref{Four1} and \eqref{Four2} leads then to the functional equation
\begin{equation}
\label{Mel_1}
P(s)\cdot Q(1-s) = 1,\ \ c_{(p,q)}< \Re(s) < C_{(p,q)}
\end{equation}
in the frequency domain.

Because the majority of elementary and special functions are particular cases of the Meijer $G$- or Fox $H$-functions, their Mellin integral transforms  can be represented in terms of the quotiens of products of the Gamma-functions. Thus, the formula \eqref{Mel_1} lets a straightforward derivation of both the Fourier kernels  and the  symmetrical Fourier kernels  in terms of the Meijer $G$- and Fox $H$-functions, see  \cite{Fox} for details.

In this paper, we focus on the Sonin kernels of the so-called general fractional integrals (GFIs) 
\begin{equation}
\label{GFI}
(\I_{(\kappa)}\, f) (x) =  (\kappa\, *\, f)(x) = \int_0^x \kappa(x-\xi)f(\xi)\, d\xi
\end{equation}
and the general fractional derivatives (GFDs) 
\begin{equation}
\label{FDR-L}
(\D_{(k)}\, f) (x) = \frac{d}{dx} (\I_{(k)}\, f) (x) = \frac{d}{dx}(k\, *\, f)(x) = \frac{d}{dx} \int_0^x k(x-\xi)f(\xi)\, d\xi  
\end{equation}
of the Laplace convolution type. 

The functions $\kappa$ and $k$ are said to form a pair of the Sonin kernels iff the function
\begin{equation}
\label{pair}
F(\xi) = (\I_{(\kappa)}\, f) (\xi),\ \xi>0
\end{equation}
solves the equation
\begin{equation}
\label{pair_1}
f(x) = (\D_{(k)}\, F) (x),\ x>0,
\end{equation}
i.e., iff the GFD \eqref{FDR-L} is a left-inverse operator to the GFI \eqref{GFI} on a suitable non-trivial space of functions (see \cite{Hil,Luc_20} for the properties that the fractional integrals and derivatives should possess). 

In its turn, the pair of the formulas  \eqref{pair} and  \eqref{pair_1} is valid iff the kernels $\kappa$ and $k$ satisfy the functional equation
\begin{equation}
\label{Son}
(\kappa \, *\, k )(x) = 1,\ x>0,
\end{equation}
where the operation $*$ stands for the Laplace convolution
\begin{equation}
\label{Lconv}
(\kappa\, *\, k)(x) = \int_0^{x}\, \kappa(x-\xi)k(\xi)\, d\xi.
\end{equation} 

In fact, the relation \eqref{Son} was introduced by Sonin in \cite{Son} as a definition of a pair of the Sonin  kernels. Nowadays it is referred to as the Sonin condition. For a given Sonin kernel $\kappa$, the kernel $k$ from the equation  \eqref{Son} is called its associated Sonin kernel. In \cite{Son}, Sonin formally proved the validity of the pair of the formulas  \eqref{pair} and  \eqref{pair_1} for the kernels $\kappa$ and $k$ that satisfy the equation  \eqref{Son}. However, an interpretation of the operators \eqref{GFI} and \eqref{FDR-L} as the GFIs and the GFDs was suggested much later by Kochubei in the paper \cite{Koch11}. In \cite{Koch11}, Kochubei also defined a special class of the Sonin kernels in terms of their Laplace integral transforms and introduced a regularized form of the GFDs:
\begin{equation}
\label{FD-C_1}
(\,_*\D_{(k)}\, f) (x) = (\D_{(k)}\, f) (x) - k(x)f(0), \ x>0.
\end{equation}

For an absolutely continuous function $f$, the right-hand side of the formula \eqref{FD-C_1} can be rewritten as
$$
(\D_{(k)}\, f) (x) - k(x)f(0) = (\I_{(k)}\, f^\prime) (x), \ x>0
$$
and thus the regularized GFD defined as in \eqref{FD-C_1} takes the form
\begin{equation}
\label{FD-C}
(\,_*\D_{(k)}\, f) (x) = (\I_{(k)}\, f^\prime) (x), \ x>0.
\end{equation}

If the Laplace integral transforms 
\begin{equation}
\label{Lap_Son}
\tilde{\kappa}(p) = \int_0^{+\infty} \kappa(x)\, e^{-px}\, dx,\ \  \tilde{k}(p) = \int_0^{+\infty} k(x)\, e^{-px}\, dx
\end{equation}
of the Sonin kernels  $\kappa$ and $k$  do exist on the domain $D_{(\kappa,k)} = \{ p\in \Com : \Re(p) > c_{(\kappa,k)}\}$, the relation \eqref{Son} is reduced to a simple formula in the Laplace domain:
\begin{equation}
\label{Son_L}
\tilde{\kappa}(p)\cdot  \tilde{k}(p)= \frac{1}{p},\  \Re(p) > c_{(\kappa,k)}.
\end{equation}

This formula can be used for derivation of some particular cases of the Sonin kernels by using the tables of the Laplace integral transforms of elementary and special functions.  However, in contrast to the case of the Fourier kernels, we cannot use the technique of the Mellin integral transform and provide a general description of the Sonin kernels in terms of the Meijer $G$- and Fox $H$-functions. In this paper, we suggest some alternative approaches for derivation of the so-called symmetrical Sonin kernels that we define in a generalized sense, see the next section for details.

In connection to the topic of this paper, it is worth mentioning that recently the GFIs and the GFDs with the Sonin kernels became a subject of active research both in Fractional Calculus (FC), see, e.g.,  \cite{Han20,Koch11,Koch19_1,LucYam16,LucYam20,Luc21a,Luc21b,Luc21c,Tar} and in its applications, see \cite{tar3}-\cite{tar7} for the models of the general  fractional dynamics,  the general non-Markovian quantum dynamics, the general non-local electrodynamics, the 
non-local classical theory of gravity, and non-local statistical mechanics, respectively, and \cite{ata,baz,gor,mis1,mis2} for the mathematical models of anomalous diffusion and linear viscoelasticity in terms of the  GFIs and the GFDs with the Sonin kernels. Thus, investigation of the general properties and particular cases of the Sonin kernels is a topic important both for the theory of FC and for its applications. 

The rest of the paper is organized as follows. In the 2nd Section, we provide some preliminary information regarding the Sonin kernels as well as their examples. The 3rd Section is devoted to derivation of the  Sonin kernels, both symmetrical and non-symmetrical,  in the time domain, whereas in the 4th Section, we employ for this aim the Laplace integral transform technique in the frequency domain.

\section{Definitions and examples}
\label{s2}

For the first time, a pair of the Sonin kernels appeared in the publications \cite{abel1} and \cite{abel2} by Abel devoted to the so-called tautochrone problem. The tautochrone is a curve for which the time taken by an object sliding without friction in uniform gravity to its lowest point is independent of its starting point on the curve. By the time of Abel's publications,  Christiaan Huygens could already prove by some advanced geometrical methods that the tautochrone is a suitable part of the cycloid. In \cite{abel1} and \cite{abel2}, Abel provided an analytical solution to a bit more general mechanical problem that is nowadays referred to as the generalized tautochrone problem: For a given function $F=F(x)$, find a curve for which the sliding time taken by an object sliding without friction  in uniform gravity depends on the position $x$ of its starting point on the curve as $F=F(x)$. In \cite{abel1} and \cite{abel2}, Abel first derived  a mathematical model of the generalized tautochrone problem in form of the following integro-differential equation (in slightly different notations)
\begin{equation}
\label{tau_1}
F(x) =  \int_0^x \frac{f^\prime(\xi)}{\sqrt{(x-\xi)}}\, d\xi
\end{equation}
and then considered a more general equation 
\begin{equation}
\label{tau_2}
F(x) =  \int_0^x \frac{f^\prime(\xi)}{(x-\xi)^{\alpha}}
\,  d\xi,  \ 0<\alpha <1.
\end{equation}

The solution to the integro-differential equation \eqref{tau_2} derived by Abel has the form (according to his tautochrone model, the condition $f(0)=0$ is valid) 
\begin{equation}
\label{tau_3}
f(x) = \frac{\sin(\alpha\, \pi)}{\pi}\, \int_0^x (x-\xi)^{\alpha-1}\, F(\xi)\, d\xi,\ x>0.
\end{equation}

It is worth mentioning that the known formula
$$
\frac{\sin(\alpha\, \pi)}{\pi} = \frac{1}{\Gamma(\alpha)\Gamma(1-\alpha)}
$$
immediately leads to a representation of the equations \eqref{tau_2}  and \eqref{tau_3}  
in terms of the so-called Caputo fractional derivative ${ }_*D_{0+}^\alpha$ and the Riemann-Liouville fractional integral $I_{0+}^\alpha$ of the order $\alpha$: 
\begin{equation}
\label{tau_4}
\frac{F(x)}{\Gamma(1-\alpha)} = ({ }_*D_{0+}^\alpha f)(x),\  \Gamma(1-\alpha)f(x) = (I_{0+}^\alpha F)(x),\ x>0,\ 0<\alpha <1.
\end{equation}

The main ingredient of the solution method  to the integro-differential equation \eqref{tau_2} invented by Abel was a simple formula for the kernels of the operators at the right-hand sides of the equations \eqref{tau_2}  and \eqref{tau_3}, i.e., of the kernels of the Riemann-Liouville fractional integral and the Riemann-Liouville or Caputo fractional derivatives  (in the modern notations)
\begin{equation}
\label{h_con}
(h_{\alpha} \, * \, h_{1-\alpha})(x) \, = \, 1,\ 0<\alpha<1,\ x>0,
\end{equation}
where the power law function $h_\alpha$ is defined as follows:
\begin{equation}
\label{h}
h_{\alpha}(x)=\frac{x^{\alpha-1}}{\Gamma(\alpha)},\ \ \alpha >0.
\end{equation}

Thus, in \cite{abel1} and \cite{abel2}, Abel derived the first and very important pair of the Sonin kernels:
\begin{equation}
\label{A_1}
\kappa(x) = h_{\alpha}(x),\ \  k(x) = h_{1-\alpha}(x),\ 0<\alpha<1,\ x>0.
\end{equation}

In particular, the tautochrone problem \eqref{tau_1} corresponds to the kernels \eqref{A_1} with $\alpha = \frac{1}{2}$. In this case, the kernels $\kappa$ and $k$ have exactly same form:
\begin{equation}
\label{A_2}
\kappa(x) = k(x) = h_{1/2}(x) = \frac{1}{\sqrt{\pi}} \frac{1}{\sqrt{x}},\ x>0.
\end{equation}

It is natural to call such kernels symmetrical Sonin kernels. However,  this is the only pair of the Sonin kernels that satisfies the relation $\kappa(x) = k(x),\ x\in \R_+$ if we suppose that the Laplace integral transforms of the kernels do exist. Indeed, in this case we get the following chain of implications from the Sonin condition \eqref{Son_L} in the frequency domain:
\begin{equation}
\label{sym_un}
\tilde{\kappa}(p)  \cdot \tilde{k}(p) = \tilde{\kappa}^2(p) = \frac{1}{p}\ \Rightarrow \ \tilde{\kappa}(p) = k(p) = \frac{1}{\sqrt{p}} \ \Rightarrow \ \kappa(x) = k(x) = h_{1/2}(x).
\end{equation}

Even the power law Sonin kernels \eqref{A_1} do not coincide unless $\alpha = 1/2$. However, they are expressed in terms of the same power law function $h_\beta$ with different values of the parameter $\beta$ and generate the Riemann-Liouville and Caputo fractional derivatives and the Riemann-Liouville fractional integral that are probably most used FC operators. 

Motivated by this example, we define the symmetrical Sonin kernels in the generalized sense:

\begin{definition} 
\label{d_1}
A pair $(\kappa,\, k)$  of the Sonin kernels is called symmetrical if both of the kernels have the same form  in terms of certain elementary or special functions with possibly different parameter values. 
\end{definition}

According to this definition,  the power law Sonin kernels \eqref{A_1} can be called symmetrical. However, not all pairs of the Sonin kernels are symmetrical. In particular, we mention here the following known pairs of non-symmetrical Sonin kernels (see \cite{Han20}, \cite{Sam}, and \cite{Zac08}):
\begin{equation}
\label{Zac}
\kappa(x) = h_{\alpha,\rho}(x),\
k(x) = h_{1-\alpha,\rho}(x) \, +\, \rho\, \int_0^x h_{1-\alpha,\rho}(\xi)\, d\xi,\  0<\alpha <1,\ \rho> 0,
\end{equation}
where
$$
h_{\alpha,\rho}(x) = \frac{x^{\alpha -1}}{\Gamma(\alpha)}\, e^{-\rho x}
$$
and
\begin{equation}
\label{Han}
\kappa(x) = h_{1-\beta+\alpha}(x)\, +\, h_{1-\beta}(x),\ 
k(x) = x^{\beta -1}\, E_{\alpha,\beta}(-x^\alpha),\ 0<\alpha < \beta <1,
\end{equation}
where $E_{\alpha,\beta}$ stands for the two-parameters Mittag-Leffler function defined by the following convergent series:
\begin{equation}
\label{ML}
E_{\alpha,\beta}(z) = \sum_{n=0}^{+\infty} \frac{z^n}{\Gamma(\alpha\, n + \beta)},\ \alpha >0,\ \beta, z\in \C.
\end{equation}

For further examples of symmetrical and non-symmetrical Sonin kernels we refer to \cite{Luc21a}, \cite{Luc21b}, \cite{Sam}, and \cite{Son}. 

It is worth mentioning that any pair $(\kappa,\, k)$ of the Sonin kernels generate a pair of the general FC operators: the 
GFI \eqref{GFI}  with the kernel $\kappa$ and the GFD  \eqref{FDR-L}  or the regularized GFD in form \eqref{FD-C_1} or \eqref{FD-C} with the kernel $k$. As an example, the Sonin kernels \eqref{Han} generate the GFI in form of a sum of two Riemann-Liouville fractional integrals
\begin{equation}
\label{GFI_2}
(\I_{(\kappa)}\, f)(x) = (I^{1-\beta+\alpha}_{0+}\, f)(x) + (I^{1-\beta}_{0+}\, f)(x) ,\ x>0
\end{equation} 
and the GFD with the Mittag-Leffler function in the kernel
\begin{equation}
\label{GFD_2}
(\D_{(k)}\, f)(x) = \frac{d}{dx}\, \int_0^x (x-\xi)^{\beta -1}\, E_{\alpha,\beta}(-(x-\xi)^\alpha)\, f(\xi)\, d\xi,\ 0<\alpha < \beta <1,\ x>0.
\end{equation} 
The regularized  GFD \eqref{FD-C} with the Mittag-Leffler function  in the kernel has then the form
\begin{equation}
\label{GFD_3}
(\,_*\D_{(k)}\, f)(x) =  \int_0^x (x-\xi)^{\beta -1}\, E_{\alpha,\beta}(-(x-\xi)^\alpha)\, f^\prime(\xi)\, d\xi,\ 0<\alpha < \beta <1,\ x>0.
\end{equation} 

On the suitable spaces of functions, the GFDs \eqref{GFD_2} and \eqref{GFD_3} are the left-inverse operators to the GFI \eqref{GFI_2} (see, e.g., \cite{Luc21a,Luc21b} for details):
\begin{equation}
\label{GFD_4}
(\D_{(k)}\, \I_{(\kappa)}\, f)(x) =  f(x), \ \ (\,_*\D_{(k)}\, \I_{(\kappa)}\, f)(x) =  f(x). 
\end{equation} 

As we see, in the case of non-symmetrical Sonin kernels, the formulas for the GFIs and the GFDs look very different. In the rest of this paper, we mainly deal with the symmetrical Sonin kernels and discuss derivation of such kernels in the time and frequency domains as well as several known and new kernels.

\section{Sonin's kernels in the time domain}
\label{s3}

In \cite{Son}, Sonin introduced an important class of the Sonin kernels in form of the products of the power law functions 
 and the analytical functions:
\begin{equation}
\label{3-3}
\kappa(x) = x^{\beta-1} \cdot \kappa_1(x),\ \kappa_1(x)=\sum_{n=0}^{+\infty}\, a_n x^n, \ a_0 \not = 0,\ 0<\beta <1,
\end{equation}
\begin{equation}
\label{3-4}
k(x) = x^{-\beta} \cdot k_1(x),\ k_1(x)=\sum_{n=0}^{+\infty}\, b_n x^n,\ 0<\beta <1.
\end{equation}
Whereas the coefficients $a_n,\ n=0,1,2,\dots$ of the analytical function $\kappa_1$ can be arbitrary chosen, the coefficients $b_n, \ n=0,1,2,\dots$ of the function $k_1$ have to satisfy the following triangular system of the linear equations
\begin{equation}
\label{3-5}
\Gamma(\beta) \Gamma(1-\beta)a_0 b_{0}=1,\ \sum_{n=0}^N\Gamma(n+\beta) \Gamma(N-n+1-\beta) a_n  b_{N-n}= 0,\ N=1,2,3\dots.
\end{equation}
The equations \eqref{3-5} are obtained by substitution of the functions $\kappa$ and $k$ defined as in \eqref{3-3} and \eqref{3-4} into the Sonin condition \eqref{Son}, interchanging the order of the integration and summation operations (that is justified by the uniform convergence of the series defining the analytical functions $\kappa_1$ and $k_1$ at any finite interval) and by applying the formula
\begin{equation}
\label{2-9}
(h_{\alpha} \, * \, h_\beta)(x) \, = \, h_{\alpha+\beta}(x),\ \alpha,\beta >0,\ x>0,
\end{equation}
that immediately follows from the well-known representation of the Euler Beta-function in terms of the Gamma-function:
$$
B(\alpha,\beta) = \int_0^1 u^{\alpha -1} (1-u)^{\beta -1}\, du = \frac{\Gamma(\alpha)\Gamma(\beta)}{\Gamma(\alpha + \beta)},\ \alpha,\beta >0.
$$

In the following lemma, a slightly more general construction of the Sonin kernels in form \eqref{3-3} and \eqref{3-4}  is provided.

\begin{lemma}
\label{l1}
Let $\kappa_1$ and $k_1$ be analytical functions and the conditions 
\begin{equation}
\label{cond_s}
0<\beta < 1,\ \alpha >0
\end{equation}
be satisfied. 

The functions 
\begin{equation}
\label{kappa}
\kappa(x) = x^{\beta-1}\cdot \kappa_1(\lambda x^\alpha),\ \kappa_1(x)=\sum_{n=0}^{+\infty}\, a_n x^n, \ a_0 \not = 0
\end{equation}
and
\begin{equation}
\label{k}
k(x) = x^{-\beta} \cdot k_1(\lambda x^\alpha),\ k_1(x)=\sum_{n=0}^{+\infty}\, b_n x^n
\end{equation}
build a pair of the Sonin kernels iff  the coefficients $a_n,\ n=0,1,2,\dots$ and $b_n, \ n=0,1,2,\dots$ of the functions $\kappa_1$ and $k_1$ satisfy the following triangular system of equations
\begin{equation}
\label{sys}
\Gamma(\beta)\Gamma(1-\beta)a_0 b_{0}=1,\ \sum_{n=0}^N\Gamma(\alpha n + \beta) \Gamma(\alpha(N-n)+1-\beta) a_n b_{N-n}= 0,\ N=1,2,3\dots.
\end{equation}
\end{lemma}

The statement formulated in Lemma \ref{l1} is proved by the same direct calculations as the ones employed for the Sonin kernels in form \eqref{3-3} and \eqref{3-4} (see, e.g., \cite{Sam}),  i.e., by substitution of the functions $\kappa$ and $k$ defined as in \eqref{kappa} and \eqref{k} into the Sonin condition \eqref{Son}, interchanging the order of the integration and summation operations and by applying the formula \eqref{2-9} for calculation of the convolution integrals. 

Now we follow another idea suggested by Sonin in his paper \cite{Son} to derive some particular cases of the symmetrical Sonin kernels in form \eqref{kappa} and \eqref{k}. First, we introduce the notations
\begin{equation}
\label{sys_1}
c_n := \Gamma(\alpha n+\beta) a_n \ \ \mbox{and} \ \ d_n:= \Gamma(\alpha n+1-\beta) b_{n}.
\end{equation}
With these notations, the equations \eqref{sys} take a simpler form:
\begin{equation}
\label{sys_2}
c_0 d_{0}=1,\ \sum_{n=0}^N  c_n d_{N-n}= 0,\ N=1,2,3\dots.
\end{equation}
Without loss of generality, in what follows, we set 
\begin{equation}
\label{null}
c_0 = d_{0}=1.
\end{equation}
Now we build a power series with the coefficients $c_n,\ n=0,1,2,\dots$ and denote its sum by $\phi = \phi(x)$:
\begin{equation}
\label{phi}
{\phi(x)} =  \sum_{n=0}^{+\infty}  c_n x^n,\ c_0 = 1.
\end{equation}
According to definition of the Cauchy product of two power law series, the sequence $d_n,\ n=0,1,2,\dots$ defined by the equations \eqref{sys_2} can be interpreted as the coefficients of the power law series of the reciprocal to the function $\phi=\phi(x)$:
\begin{equation}
\label{1/phi}
\frac{1}{\phi(x)} =  \sum_{n=0}^{+\infty}  d_n x^n,\ d_0 = 1.
\end{equation}

Summarising the arguments presented above, we thus arrive at a method for construction of particular cases of the Sonin kernels by means of the known coefficients of the Taylor series for some analytical functions and their reciprocals. In particular, this method leads to the symmetrical Sonin kernels in the case when the functions $\phi(x)$ and $\frac{1}{\phi(x)}$ defined by the series \eqref{phi} and \eqref{1/phi}, respectively, has the  same form  in terms of certain elementary or special functions with possibly different parameter values. 

\begin{remark}
\label{rem1}
It is worth mentioning that  under the conditions \eqref{sys}, the functions in form \eqref{3-3} and \eqref{3-4} are still the Sonin kernels  if the functions $\kappa_1$ and $k_1$ are represented by the power series with the finite convergence radii greater than zero, i.e., if they are not  analytical functions anymore (see \cite{Wick} for the proofs). 
 The arguments presented in \cite{Wick} remain valid also for the kernels \eqref{kappa} and \eqref{k} from Lemma \ref{l1} and thus we can relax the condition of analyticity of the functions $\kappa_1$ and $k_1$ from the right-hand sides of the formulas \eqref{kappa} and \eqref{k}, respectively, to the condition that the convergence radii of the Taylor  series  of the functions $\kappa_1$ and $k_1$ are greater than zero.
\end{remark}

In what follows, we apply the method presented above for derivation of some Sonin kernels, both already known and the new ones. 



\begin{example}
    \label{ex_3_1}

In the formula \eqref{phi}, we set $\phi(x) = \exp(x)$. Then
\begin{equation}
\label{phi_1}
{\phi(x)} =  \exp(x)=  \sum_{n=0}^{+\infty}  \frac{1}{n!}\, x^n,\ \ 
\frac{1}{\phi(x)} =\exp(-x) =  \sum_{n=0}^{+\infty} \frac{(-1)^n}{n!}\, x^n.
\end{equation}

Thus, 
\begin{equation}
\label{cn_1}
c_n = \frac{1}{n!},\ \ d_n = \frac{(-1)^n}{n!},\ n=0,1,2,\dots
\end{equation}
and (see the relations \eqref{sys_1})
\begin{equation}
\label{an_1}
a_n = \frac{c_n}{\Gamma(\alpha n + \beta)} = \frac{1}{n!\Gamma(\alpha n + \beta) }, \ n=0,1,2,\dots ,
\end{equation}
\begin{equation}
\label{bn_1}
b_n =\frac{d_n}{\Gamma(\alpha n +1-\beta)} 
= \frac{(-1)^n}{n!\Gamma(\alpha n +1-\beta)},\ n=0,1,2,\dots.
\end{equation}

The Sonin kernels \eqref{kappa} and \eqref{k} from Lemma \ref{l1} with the coefficients as in \eqref{an_1} and \eqref{bn_1}, respectively, take then the form ($\beta \in (0,\, 1),\ \alpha >0$)
\begin{equation}
\label{kappa_1}
\kappa(x) = x^{\beta-1}\, \sum_{n=0}^{+\infty} \frac{(\lambda x^\alpha)^n}{n!\Gamma(\alpha n + \beta) } = 
x^{\beta-1}\, W_{\alpha,\beta}(\lambda x^\alpha)
\end{equation}
and
\begin{equation}
\label{k_1}
k(x) = x^{-\beta}\, \sum_{n=0}^{+\infty} \frac{(-\lambda x^\alpha)^n}{n!\Gamma(\alpha n + 1-\beta) } = x^{-\beta}\, W_{\alpha,1-\beta}(-\lambda x^\alpha),
\end{equation}
where the Wright function $W_{\alpha,\beta}(z)$ is defined by the convergent series (see \cite{Luc_19} for its properties and applications)
\begin{equation}
\label{wright}
W_{\alpha,\beta}(z) := \sum_{n=0}^{+\infty} \frac{z^n}{n!\Gamma(\alpha n+\beta) },\ z,\beta  \in \C,\ \alpha >-1.
\end{equation}
\end{example}

To the best knowledge of the author, the pair \eqref{kappa_1}, \eqref{k_1} of the symmetrical (in the sense of Definition \ref{d_1}) Sonin kernels has not yet been reported in the literature. However, in his paper \cite{Son}, Sonin mentioned a particular case of the kernels \eqref{kappa_1}, \eqref{k_1} with $\alpha =1$. In this case, the Sonin kernels can be expressed in terms of the  cylinder or Bessel functions, see \cite{Sam,Son} for details.


\begin{example}
\label{ex_3_2}    

In this example, we derive another symmetrical pair of the Sonin kernels generated by the function $\phi(x) = (1+x)^{-\gamma}$. The Taylor series for 
$\phi(x)$ and its reciprocal are well-known: 
\begin{equation}
\label{phi_2}
{\phi(x)} =  (1+x)^{-\gamma} =  \sum_{n=0}^{+\infty}  \frac{(-1)^n (\gamma)_n}{n!}\, x^n,\ \ 
\frac{1}{\phi(x)} =(1+x)^{\gamma} =  \sum_{n=0}^{+\infty} \frac{(-1)^n (-\gamma)_n}{n!}\, x^n,
\end{equation}
where the Pochhammer symbol $(\gamma)_n$ is defined as follows:
$$
(\gamma)_n := \frac{\Gamma(\gamma +n)}{\Gamma(\gamma)} = \gamma \cdot (\gamma +1)\cdot \dots \cdot (\gamma + n-1). 
$$

As we see, the functions  $\phi(x) = (1+x)^{-\gamma}$ and $\frac{1}{\phi(x)} =(1+x)^{\gamma} $ generate the coefficients
\begin{equation}
\label{cn_2}
c_n = \frac{(-1)^n (\gamma)_n}{n!},\ \ d_n = \frac{(-1)^n (-\gamma)_n}{n!},\ n=0,1,2,\dots ,
\end{equation}
and 
\begin{equation}
\label{an_2}
a_n = \frac{c_n}{\Gamma(\alpha n+\beta)} = \frac{(-1)^n (\gamma)_n}{n!\Gamma(\alpha n+\beta) }, \ n=0,1,2,\dots ,
\end{equation}
\begin{equation}
\label{bn_2}
b_n =\frac{d_n}{\Gamma(\alpha n+1-\beta)} 
= \frac{(-1)^n (-\gamma)_n}{n!\Gamma(\alpha n+1-\beta)},\ n=0,1,2,\dots.
\end{equation}

The Sonin kernels \eqref{kappa} and \eqref{k}  take then the form ($\beta \in (0,\, 1),\ \alpha >0$)
\begin{equation}
\label{kappa_2}
\kappa(x) = x^{\beta-1}\, \sum_{n=0}^{+\infty} \frac{(-1)^n (\gamma)_n}{n!\Gamma(\alpha n+\beta) } (\lambda x^\alpha)^n  = 
x^{\beta-1}\, E_{\alpha,\beta}^\gamma(-\lambda x^\alpha)
\end{equation}
and
\begin{equation}
\label{k_2}
k(x) = x^{-\beta}\, \sum_{n=0}^{+\infty} \frac{(-1)^n (-\gamma)_n}{n!\Gamma(\alpha n + 1-\beta) } (\lambda x^\alpha)^n  = 
x^{-\beta}\, E_{\alpha,1-\beta}^{-\gamma}(-\lambda x^\alpha),
\end{equation}
where the three-parameters Mittag-Leffler function or the Prabhakar function $E_{\alpha,\beta}^\gamma(z)$ is defined by the convergent series (see, e.g.,  \cite{Giu,Prab_2} for its properties and applications)
\begin{equation}
\label{Prab}
E_{\alpha,\beta}^\gamma(z) := \sum_{n=0}^{+\infty} \frac{(\gamma)_n}{n!\Gamma(\alpha n+\beta) },\ z,\beta,\gamma  \in \C,\ \alpha >0.
\end{equation}
\end{example}

Please note that the function \eqref{Prab} can be also interpreted as a particular case of the so-called generalized Wright or Fox-Wright function, see \cite{Luc_19} for details. 

In \cite{Son}, Sonin derived a particular case of the kernels \eqref{kappa_2}, \eqref{k_2} with $\alpha =1$ in form of the power law series. As mentioned in \cite{Sam}, the kernels \eqref{kappa_2}, \eqref{k_2} with $\alpha =1$ can be represented in terms of the confluent hypergeometric function or the Kummer function  ${ }_1 F_1$:
\begin{equation}
\label{kappa_2_1}
\kappa(x) = x^{\beta-1}\, \sum_{n=0}^{+\infty} \frac{(-1)^n (\gamma)_n}{n!\Gamma(n+\beta) } (\lambda x)^n  = 
\frac{x^{\beta-1}}{\Gamma(\beta)}\, { }_1 F_1(\gamma;\beta;-\lambda x),
\end{equation}
\begin{equation}
\label{k_2_1}
k(x) = x^{-\beta}\, \sum_{n=0}^{+\infty} \frac{(-1)^n (-\gamma)_n}{n!\Gamma( n + 1-\beta) } (\lambda x)^n  = 
\frac{x^{-\beta}}{\Gamma(1-\beta)}\, { }_1 F_1(-\gamma;1-\beta;-\lambda x),
\end{equation}
where the Kummer function ${ }_1 F_1$ is defined by the convergent series
\begin{equation}
\label{Kum}
{ }_1 F_1(\gamma;\beta;z) := \sum_{n=0}^{+\infty} \frac{(\gamma)_n}{(\beta)_n }\frac{z^n}{n!},\ z,\beta,\gamma  \in \C.
\end{equation}

The integral operator in form \eqref{GFI} (that we called the GFI with the kernel $\kappa$) and the integro-differential operator \eqref{FDR-L} (that we called the GFD with the kernel $k$) with the Sonin kernels \eqref{kappa_2_1} and \eqref{k_2_1} were considered by Prabhakar  in his paper \cite{Prab_1} (see also \cite{Samko}, formulas (37.1) and (37.31)). In particular, in \cite{Prab_1}, the GFD \eqref{FDR-L} with the kernel \eqref{k_2_1} was shown to be a left-inverse operator to the GFI \eqref{FDR-L} with the kernel \eqref{kappa_2_1}.

Two years after publication of his paper \cite{Prab_1}, in \cite{Prab_2}, Prabhakar also studied the general case of the GFI \eqref{GFI} and the GFD \eqref{FDR-L} with the Sonin kernels  \eqref{kappa_2} and \eqref{k_2}, respectively, in terms of the three-parameters Mittag-Leffler function. For this reason, nowadays this function is often referred to as the Prabhakar function. For a detailed presentation of the theory and applications of these operators called nowadays the Prabhakar fractional integral and derivative, we refer the interested reader to the recent survey \cite{Giu}. 

As mentioned in Introduction, investigation of the GFIs and the GFDs with the Sonin kernels is a very recent topic initialized by Kochubei in the paper \cite{Koch11} published in 2011. In the earlier publications by 
 Prabhakar and other authors who investigated the operators of type  \eqref{GFI} and  \eqref{FDR-L} with the kernels in terms of the Bessel function, Kummer function, three-parameters Mittag-Leffler function, etc., these operators were introduced and studied independently of the theory of the GFIs and the GFDs with the Sonin kernels. However, recently one started to interpret these earlier results from the viewpoint of the GFIs and the GFDs, see, e.g., \cite{Giu_1} for discussion of a connection between the General Fractional Calculus (GFC) and  the Prabhakar fractional integral and derivative. 


\begin{example}
\label{ex_3_3}
In this example, we consider a pair of the Sonin kernels generated by the function $\phi(x) = \exp(x)(1+x)^{-\gamma}$. The Taylor series for the function
$\phi(x)$ and its reciprocal can be obtained by multiplication of the Taylor series from Examples \ref{ex_3_1} and \ref{ex_3_2}:
\begin{equation}
\label{phi_3}
{\phi(x)} = \exp(x) (1+x)^{-\gamma} = \left( \sum_{n=0}^{+\infty}  \frac{1}{n!}\, x^n\right) \cdot \left( \sum_{n=0}^{+\infty}  \frac{(-1)^n (\gamma)_n}{n!}\, x^n\right) = \sum_{n=0}^{+\infty} c_{n,\gamma}\, x^n,
\end{equation}
where 
\begin{equation}
\label{cng}
c_{n,\gamma} = \sum_{m=0}^n \frac{(-1)^m(\gamma)_m}{m! (n-m)!},
\end{equation}
\begin{equation}
\label{phi_3_1}
\frac{1}{\phi(x)} = \exp(-x) (1+x)^{\gamma} = \left( \sum_{n=0}^{+\infty}  \frac{(-1)^n}{n!}\, x^n\right) \cdot \left( \sum_{n=0}^{+\infty}  \frac{(-1)^n (-\gamma)_n}{n!}\, x^n\right) = \sum_{n=0}^{+\infty} d_{n,\gamma}\, (-x)^n,
\end{equation}
where 
\begin{equation}
\label{dng}
d_{n,\gamma} = \sum_{m=0}^n \frac{(-\gamma)_m}{m! (n-m)!}.
\end{equation}

Taking into account the relations \eqref{sys_1} and by applying Lemma \ref{l1}, we arrive at the following new pair of the Sonin kernels ($\beta \in (0,\, 1),\ \alpha >0$):
\begin{equation}
\label{kappa_3}
\kappa(x) = x^{\beta-1}\, \sum_{n=0}^{+\infty} \frac{c_{n,\gamma}}{\Gamma(\alpha n+\beta) } (\lambda x^\alpha)^n,\ 
k(x) = x^{-\beta}\, \sum_{n=0}^{+\infty} \frac{d_{n,\gamma}}{\Gamma(\alpha n + 1-\beta) } (-\lambda x^\alpha)^n,
\end{equation}
where the the coefficients $c_{n,\gamma}$ and $d_{n,\gamma}$ are given by the formulas \eqref{cng} and \eqref{dng}, respectively. 
\end{example}

In the rest of this section, we discuss an important  extension of the Sonin kernels \eqref{3-3} and \eqref{3-4} in form of the products of the power law functions 
and the analytical functions to  the so-called convolution series. The convolution series are a far reaching generalization of the power law series that was recently introduced in the framework of the general FC (see, e.g., \cite{Luc21b,Luc22}). First, we remind the readers on their definition and properties that we need for the further discussions. 

Let a function $f$ belong to the space  $C_{-1}(0,+\infty)\, = \, \{\phi:\ \phi(x) = x^{p}\phi_1(x),\ x>0,\ p>-1,\ \phi_1\in C[0,+\infty)\}$  and be represented in the form
\begin{equation}
\label{rep}
f(x) = x^{p}f_1(x),\ x>0,\ p>0,\ f_1\in C[0,+\infty)
\end{equation} 
and let the convergence radius of the power  series
\begin{equation}
\label{ser}
\Sigma(z) = \sum^{+\infty }_{n=0}a_{n}\, z^n,\ a_{n}, z\in \C
\end{equation}
be non-zero.

The convolution series generated by the function $f$ has the form
\begin{equation}
\label{conser}
\Sigma_f(x) = \sum^{+\infty }_{n=0}a_{n}\, f^{<n+1>}(x),
\end{equation}
where $f^{<n>}$ stands for the convolution powers
\begin{equation}
\label{I-6}
f^{<n>}(x):= \begin{cases}
1,& n=0,\\
f(x), & n=1,\\
(\underbrace{f* \ldots * f}_{n\ \mbox{times}})(x),& n=2,3,\dots .
\end{cases}
\end{equation}

As was shown in \cite{Luc21b,Luc22}, the convolution series \eqref{conser}
is convergent for all $x>0$ and defines a function  from the space $C_{-1}(0,+\infty)$. 
Moreover, the series 
\begin{equation}
\label{conser_p}
x^{1-\alpha}\, \Sigma_f(x) = \sum^{+\infty }_{n=0}a_{n}\, x^{1-\alpha}\, f^{<n+1>}(x),\ \ \alpha = \min\{p,\, 1\}
\end{equation}
is uniformly convergent for $x\in [0,\, X]$ for any $X>0$.

In particular, any power series can be represented in form of the convolution series \eqref{conser} generated by the function $f(x) = h_1(x) = \frac{x^0}{\Gamma(1)} \equiv 1,\ x\ge 0$:
\begin{equation}
\label{I-ps}
\Sigma_{h_1}(x) = \sum_{n=0}^{+\infty} a_n\, h_1^{<n+1>}(x) = \sum_{n=0}^{+\infty} a_n\, h_{n+1}(x) = \sum^{+\infty }_{n=0}a_{n}\, \frac{x^n}{n!}.
\end{equation}

For examples and applications of the convolution series in form \eqref{conser} we refer to \cite{Luc21b,Luc22}. In the following theorem, we present a construction of the Sonin kernels in terms of the convolution series.

\begin{theorem}
\label{t1}
Let $(\kappa_1,\ k_1)$ be a pair of the Sonin kernels and $f$ be any function from the space $C_{-1}(0,+\infty)$ that is not identically equal to zero and let the convergence radii of the power  series
\begin{equation}
\label{ser_1}
\Sigma_a(z) = \sum^{+\infty }_{n=0}a_{n}\, z^n,\ a_{n}, z\in \C \ \ \mbox{and} \ \ \Sigma_b(z) = \sum^{+\infty }_{n=0}b_{n}\, z^n,\ b_{n}, z\in \C
\end{equation}
be non-zero. 

Then the functions
\begin{equation}
\label{3-3_1}
\kappa (x) = a_0\kappa_1(x) + (\kappa_1\, *\, \sum_{n=1}^{+\infty}\, a_n\, f^{<n>})(x), \ a_0 \not = 0,
\end{equation}
\begin{equation}
\label{3-4_1}
k(x) = b_0 k_1(x) + (k_1\, *\, \sum_{n=1}^{+\infty}\, b_n\, f^{<n>})(x)
\end{equation}
build a pair of the Sonin kernels iff the conditions 
\begin{equation}
\label{condab}
a_0 b_{0}=1,\ \sum_{m=0}^n  a_m b_{n-m}= 0,\ n=1,2,3\dots 
\end{equation}
hold true.
\end{theorem}

\begin{proof}
By the Sonin condition \eqref{Son} and by definition of the convolution powers we immediately get the relations 
$$
(a_0\kappa_1\, *\, b_0k_1)(x) = a_0b_0,\    (a_0\kappa_1\, *\, (b_m k_1 \, *\, f^{<m>}))(x) =a_0 b_m(1 \, *\, f^{<m>})(x),
$$
$$
 ((\kappa_1 \, *\, a_i f^{<i>})\, *\, b_0 k_1)(x) = a_i b_0(1 \, *\, f^{<i>})(x),\  
$$
$$
  ((\kappa_1 \, *\, a_i f^{<i>})\, *\, (k_1 \, *\, b_m f^{<m>}))(x) = 
a_i b_m (1 \, *\, f^{<i+m>})(x).
$$
To calculate the convolution $(\kappa\, *\, k)(x)$ of the functions $\kappa$ and $k$ defined by \eqref{3-3_1} and \eqref{3-4_1}, respectively, we interchange the orders of integration and summation (that is allowed because of the uniform convergence of the convolution series \eqref{conser_p}) and use the relations from above. Thus we arrive at the representation
$$
(\kappa\, *\, k)(x) = a_0b_0 +\sum_{n=1}^{+\infty} \left(\sum_{m=0}^n a_mb_{n-m}\right) (1\, *\, f^{<n>})(x).
$$ 
This formula  ensures that the pair of the functions $\kappa$ and $k$ given by the relations \eqref{3-3_1} and 
\eqref{3-4_1} are the Sonin kernels iff the  conditions \eqref{condab} are satisfied. 
\end{proof}

\begin{remark}
\label{r1}
The Sonin kernels \eqref{3-3} and \eqref{3-4} in form of the products of the power law functions $x^{\beta-1}$ and $x^{-\beta}$ with $0<\beta <1$ and the analytical functions are a particular case of the general convolution kernels \eqref{3-3_1} and \eqref{3-4_1} with $\kappa_1(x) = h_{\beta}(x)$, $k_1(x)= h_{1-\beta}(x)$ and $f(x) = h_1(x) \equiv 1$. Other examples can be constructed by employing different Sonin kernels $\kappa_1,\ k_1$ and the function $f$ that generates the convolution series in the formulas \eqref{3-3_1} and \eqref{3-4_1}. In particular, the Sonin kernels \eqref{kappa} and \eqref{k} presented in Lemma \ref{l1} correspond  to the convolution kernels \eqref{3-3_1} and \eqref{3-4_1} with $\kappa_1(x) = h_{\beta}(x)$, $k_1(x) = h_{1-\beta}(x)$ ($0<\beta <1$) and $f(x) = \lambda\, h_\alpha(x),\ \alpha >0$.
\end{remark}

\section{Sonin's kernels in the Laplace domain}
\label{s4}

As already mentioned in Introduction, the Laplace integral transforms of any pair of the Sonin kernels (if they exist) are connected each to other by the relation
\eqref{Son_L} 
that can be easily solved for, say, $\tilde{k}$:
\begin{equation}
\label{Son_L_1}
\tilde{k}(p)= \frac{1}{p\tilde{\kappa}(p)},\  \Re(p) > c_{(\kappa,k)}.
\end{equation}

The last relation means that two Sonin kernels, $\kappa$ and $k$, are symmetrical iff the functions $\tilde{\kappa}(p)$ and $\frac{1}{p\tilde{\kappa}(p)}$ can be represented in terms of the same elementary or special function with possibly different parameters values. In what follows, we consider some examples of the known and new symmetrical Sonin kernels that possess this property. 


\begin{example}
\label{ex_4_1}
We start with the well-know case of the Laplace transform $\tilde{\kappa}(p)$ in terms of a power function
\begin{equation}
\label{Lap}
\tilde{\kappa}(p) = p^{-\alpha},\  \alpha>0,\ \Re(p) >0.
\end{equation}
Then 
\begin{equation}
\label{Lap_k}
\tilde{k}(p)= \frac{1}{p\tilde{\kappa}(p)} = p^{\alpha-1},\  \alpha<1, \  \Re(p) >0.
\end{equation}

The Laplace transforms $\tilde{\kappa}(p)$ and $\tilde{k}(p)$ are expressed in terms of the power function with different exponents. The basic formula
\begin{equation}
\label{power_L}
\tilde{h_\gamma}(p) = p^{-\gamma},\  \gamma >0, \ \Re(p) > 0
\end{equation}
for the Laplace transform of the power function $h_\gamma(t) = t^{\gamma -1}/\Gamma(\gamma)$ leads then to the well-known pair of the symmetrical Sonin kernels in the time domain:
$$
\kappa(t) = h_{\alpha}(t),\ \ k(t) = h_{1-\alpha}(t),\ 0<\alpha<1.
$$
\end{example}


\begin{example}
\label{ex_4_2}
Now let us consider the case
\begin{equation}
\label{Lap_2}
\tilde{\kappa}(p) = p^{-\beta}\exp(\lambda\, p^{-\alpha}), \ \alpha>0,\ \beta >0,\ \Re(p) > 0.
\end{equation}
The Laplace transform of the associated kernel $k$ has then a similar form
\begin{equation}
\label{Lap_2_1}
\tilde{k}(p) = \frac{1}{p\kappa(p)}= p^{\beta-1}\exp(-\lambda\, p^{-\alpha}), \ \alpha >0,\ \beta <1, \ \Re(p) > 0.
\end{equation}

We note here that the functions \eqref{Lap_2} and \eqref{Lap_2_1}  with $\lambda = 0$ have been already considered in Example \ref{ex_4_2}.

To represent the function $\kappa$ in the time domain, we employ the series representation of its Laplace transform
\begin{equation}
\label{Lap_2_2}
\tilde{\kappa}(p) =p^{-\beta}\exp(\lambda\, p^{-\alpha})=  \sum_{n=0}^{+\infty}\frac{\lambda^n}{n!}p^{- \alpha n -\beta}
\end{equation}
and the formula \eqref{power_L} and thus arrive to the expression
\begin{equation}
\label{time_2}
\kappa(x)  = \sum_{n=0}^{+\infty} \frac{ \lambda^n}{n!} h_{\alpha n+\beta}(x)  = \sum_{n=0}^{+\infty} \frac{ \lambda^n}{n!} \frac{x^{\alpha n+\beta-1}}{\Gamma(\alpha n +\beta)} = x^{\beta -1} W_{\alpha,\beta}(\lambda x^\alpha), \ \alpha>0,\ \beta >0,
\end{equation}
where $W_{\alpha,\beta}$ is the Wright function \eqref{wright}. Using the same method for the Laplace transform 
$\tilde{k}(p)$ given by \eqref{Lap_2_1}, we get a representation of the associated kernel $k$ in the time domain:
\begin{equation}
\label{time_2_1}
k(x)  = x^{-\beta}  W_{\alpha,1-\beta}(-\lambda x^\alpha),\ \alpha >0,\ \beta <1, 
\end{equation}
where $W_{\alpha,\beta}$ is the Wright function \eqref{wright}. 
\end{example}

To the best knowledge of the author, the pair of the Sonin kernels in form \eqref{time_2}, \eqref{time_2_1} was not yet reported in the literature. However, we already derived the same Sonin kernels in Section \ref{s3}, see Example \ref{ex_3_1}. 


\begin{example}
\label{ex_4_3}
In this example, we set 
\begin{equation}
\label{Lap_3}
\tilde{\kappa}(p) = 
\frac{p^{\alpha \gamma -\beta}}{(p^\alpha +\lambda)^\gamma},\ \alpha>0,\ \beta >0,\  |\lambda p^{-\alpha}|<1. 
\end{equation}
The Laplace transform of the associated kernel $k$ has the same form with different parameters values:
\begin{equation}
\label{Lap_3_1}
\tilde{k}(p) = \frac{1}{p\kappa(p)}= \frac{(p^\alpha +\lambda)^\gamma }{p^{\alpha \gamma -\beta +1}} = \frac{p^{-\alpha \gamma +\beta -1}}{(p^\alpha +\lambda)^{-\gamma} }, \ \beta < 1,\ \alpha>0,\ |\lambda p^{-\alpha}|<1. 
\end{equation}

For $\lambda= 0$, the kernels \eqref{Lap_3} and \eqref{Lap_3_1} are reduced to the kernels that have been discussed in Example \ref{ex_4_1}.

Using the power series representation
\begin{equation}
\label{Lap_3+}
\tilde{\kappa}(p) = \frac{p^{\alpha \gamma -\beta}}{(p^\alpha +\lambda)^\gamma} = p^{-\beta}( 1 + \lambda p^{-\alpha})^{-\gamma} = \sum_{n=0}^{+\infty}\frac{(-1)^n (\gamma)_n \lambda^n}{n!}p^{-\alpha n -\beta}
\end{equation}
and the formula \eqref{power_L}, we get the following result in the time domain:
\begin{equation}
\label{time_3}
\kappa(x)  = \sum_{n=0}^{+\infty} \frac{(-1)^n (\gamma)_n  \lambda^n }{n! }h_{\alpha n +\beta}(x) =  x^{\beta -1} E_{\alpha,\beta}^\gamma(-\lambda x^\alpha),\ \alpha >0,\ \beta >0,
\end{equation}
where $E_{\alpha,\beta}^\gamma$ is the three-parameters Mittag-Leffler function defined as in \eqref{Prab}. 

The same method applied to the Laplace transform \eqref{Lap_3_1} easily leads to the formula
\begin{equation}
\label{time_3_1}
k(x)  = x^{-\beta} E_{\alpha,1-\beta}^{-\gamma}(-\lambda x^\alpha),\ \alpha >0,\ \beta < 1,
\end{equation}
where $E_{\alpha,\beta}^\gamma$ is the three-parameters Mittag-Leffler function defined as in \eqref{Prab}. 
\end{example}

The pair of the symmetrical Sonin kernels \eqref{time_3} and \eqref{time_3_1} has been already discussed in Section \ref{s3}, see Example \ref{ex_3_2}. 


\begin{example}
\label{ex_4_4}
In this example, we consider the Sonin kernel $\kappa$ with the Laplace integral transform $\tilde{\kappa}$ in the form
\begin{equation}
\label{Lap_4}
\tilde{\kappa}(p) = p^{-\beta} \frac{p^{\alpha_1 \gamma}}{(p^{\alpha_1} +\lambda_1)^\gamma} \exp(\lambda_2\, p^{-\alpha_2}) = p^{\alpha_1 \gamma -\beta}(p^{\alpha_1} +\lambda_1)^{-\gamma}\exp(\lambda_2\, p^{-\alpha_2}) ,
\end{equation}
where the parameters and the Laplace variable satisfy the conditions $\alpha_1>0,\ \alpha_2>0,\ \beta >0,\ \Re(p) > 0,\ |\lambda_1 p^{-\alpha_1}|<1$. 

Please note that we already considered two important particular cases of the kernels in form \eqref{Lap_4}, namely, the case $\lambda_1 = 0$ in Example \ref{ex_4_2} and the case $\lambda_2 = 0$ in Example \ref{ex_4_3}.

The Laplace transform of the associated kernel $k$ has then a similar form
\begin{equation}
\label{Lap_4_1}
\tilde{k}(p) = \frac{1}{p\kappa(p)}=  p^{-\alpha_1 \gamma +\beta-1}  (p^{\alpha_1} +\lambda_1)^{\gamma} \exp(-\lambda_2\, p^{-\alpha_2}),
\end{equation}
where the parameters and the Laplace variable satisfy the conditions $\alpha_1>0,\ \alpha_2>0,\  \beta <1,\ \Re(p) > 0,\ |\lambda_1 p^{-\alpha_1}|<1$.

The series representations of the functions \eqref{Lap_4} and \eqref{Lap_4_1} immediately follow from the formulas \eqref{Lap_2_2} and \eqref{Lap_3+}:
\begin{equation}
\label{Lap_4_2}
\tilde{\kappa}(p) = \sum_{m,n=0}^{+\infty}\frac{(-\lambda_1)^m (\gamma)_m \lambda_2^n }{m! n!}p^{-\alpha_1 m- \alpha_2 n-\beta},
\end{equation}
\begin{equation}
\label{Lap_4_3}
\tilde{k}(p) = \frac{1}{p\kappa(p)} =  \sum_{m,n=0}^{+\infty}\frac{(-\lambda_1)^m (-\gamma)_m (-\lambda_2)^n }{m! n!}p^{-\alpha_1 m - \alpha_2 n+\beta-1}.
\end{equation}
The last two formulas are valid under the conditions $\alpha_1>0,\ \alpha_2>0,\ \ 0<\beta <1,\  |\lambda_1 p^{-\alpha_1}|<1,\ \Re(p) > 0$. 

Now we apply the formula \eqref{power_L} and get the following representations of the Sonin kernels $\kappa$ and $k$ in the time domain:
\begin{equation}
\label{time_4_4}
\kappa(x)  = \sum_{m,n=0}^{+\infty}\frac{ (-\lambda_1)^m (\gamma)_m \lambda_2^n }{m! n!}h_{\alpha_1 m + \alpha_2 n+\beta}(x) = x^{\beta-1}\, \phi_3(\gamma;(\alpha_1,\alpha_2;\beta);-\lambda_1x^{\alpha_1},\lambda_2x^{\alpha_2}),
\end{equation}
$$
k(x)  = \sum_{m,n=0}^{+\infty}\frac{(-\lambda_1)^m (-\gamma)_m (-\lambda_2)^n }{m! n!}h_{\alpha_1 m + \alpha_2 n -\beta+1}(x) = 
$$
\begin{equation}
\label{time_4_5}
x^{-\beta}\, \phi_3(-\gamma;(\alpha_1,\alpha_2;1-\beta);-\lambda_1x^{\alpha_1},-\lambda_2x^{\alpha_2}),
\end{equation}
where $\alpha_1>0,\ \alpha_2>0,\ 0<\beta <1$ and $\phi_3$ is a new special function of the hypergeometric type in two variables defined by the convergent series
\begin{equation}
\label{W_3}
\phi_3(\gamma;(\alpha_1,\alpha_2;\beta);y,z):= \sum_{m,n=0}^{+\infty} \frac{(\gamma)_m}{\Gamma(\alpha_1 m + \alpha_2 n+\beta)}\, \frac{y^mz^n}{m! n!},\ \alpha_1>0,\ \alpha_2>0,\ \gamma,\beta,y,z\in \Com.
\end{equation}
\end{example}

The denotation $\phi_3$ is motivated by a particular case of this function for $\alpha_1 =\alpha_2 =1$ that is reduced to the known  Horn function $\Phi_3$:
\begin{equation}
\label{Horn_r}
\phi_3(\gamma;(1,1;\beta);y,z)= \sum_{m,n=0}^{+\infty} \frac{(\gamma)_m}{\Gamma((m+n)+\beta)}\, \frac{y^mz^n}{m! n!} = \frac{1}{\Gamma(\beta)} 
\Phi_3(\gamma;\beta;y,z),
\end{equation}
where $\Phi_3$ is one of the Horn functions defined by the double confluent series of the hypergeometric type (see the formula (22) in section 5.7.1 in \cite{Bat}):
\begin{equation}
\label{Horn}
\Phi_3(\gamma;\beta;y,z):= \sum_{m,n=0}^{+\infty} \frac{(\gamma)_m}{(\beta)_{m+n}}\, \frac{y^mz^n}{m! n!}.
\end{equation}

Thus,  for $\alpha_1 = \alpha_2 = 1$, the pair \eqref{time_4_4} and \eqref{time_4_5} of the symmetrical Sonin kernels  is expressed in terms of the Horn function $\Phi_3$:
\begin{equation}
\label{time_4_6}
\kappa(x)  =  h_\beta(x)\, \Phi_3(\gamma;\beta;-\lambda_1x,\lambda_2x), \  \beta >0,
\end{equation}
\begin{equation}
\label{time_4_7}
k(x)  =  h_{1-\beta}(x)\, \Phi_3(-\gamma;1-\beta;-\lambda_1x, -\lambda_2x), \ \beta <1.
\end{equation}

The representations \eqref{time_4_6} and \eqref{time_4_7} can be also obtained by employing the formula
2.2.3.16 from \cite{Pru} that is valid under the conditions $\Re(\mu+\nu)<0, \Re(p)>\max\{0,-\Re(b)\}$:
\begin{equation}
\label{inv}
\{ \mathcal{L}^{-1} p^{\mu}(p+b)^\nu \exp(a/p)\}(x) = h_{-\mu-\nu}(x)\, \Phi_3(-\nu;-\mu-\nu;-bx,ax),
\end{equation}
where  $\Phi_3$ is the Horn function defined by \eqref{Horn} and   $\{ \mathcal{L}^{-1}\, f(p)\}(x)$ stands for the inverse Laplace transform of the function $f$ at the point $x>0$.

To the best knowledge of the author, both the functions \eqref{time_4_4} and \eqref{time_4_5} and the functions \eqref{time_4_6} and \eqref{time_4_7} are  new symmetrical Sonin kernels not yet mentioned in the literature.


\begin{example}
\label{ex_4_5}
In this last example,  the Laplace integral transform $\tilde{\kappa}$ of the Sonin kernel $\kappa$ is a product of three different power law functions:
\begin{equation}
\label{Lap_5}
\tilde{\kappa}(p) = p^{-\beta}  \frac{p^{\alpha_1 \gamma_1}}{(p^{\alpha_1} +\lambda_1)^{\gamma_1}} \frac{p^{\alpha_2 \gamma_2}}{(p^{\alpha_2} +\lambda_2)^{\gamma_2}} = p^{\alpha_1 \gamma_1+\alpha_2\gamma_2 -\beta} (p^{\alpha_1} +\lambda_1)^{-\gamma_1} (p^{\alpha_2} +\lambda_2)^{-\gamma_2} ,
\end{equation}
where the parameters and the Laplace variable satisfy the conditions $\alpha_1>0,\ \alpha_2>0,\ \beta >0,\ \Re(p) > 0,\ |\lambda_1 p^{-\alpha_1}|<1,\ |\lambda_2 p^{-\alpha_2}|<1$.

In the case $\lambda_1 = 0$ or $\lambda_2 =0$, the kernel \eqref{Lap_5} is reduced to the kernel already considered in Example \ref{ex_4_3}. 

The Laplace transform of the associated kernel $k$ takes a similar form
\begin{equation}
\label{Lap_5_1}
\tilde{k}(p) = \frac{1}{p\kappa(p)}=  p^{-\alpha_1 \gamma_1-\alpha_2 \gamma_2 +\beta-1} (p^{\alpha_1} +\lambda_1)^{\gamma_1} (p^{\alpha_2} +\lambda_2)^{\gamma_2},
\end{equation}
where the parameters and the Laplace variable satisfy the conditions $\alpha_1>0,\ \alpha_2>0,\ \beta <1,\ \Re(p) > 0,\ |\lambda_1 p^{-\alpha_1}|<1,\ |\lambda_2 p^{-\alpha_2}|<1$.

Applying the formula \eqref{Lap_3+}, we get the following series representations of the functions as in the formulas \eqref{Lap_5} and \eqref{Lap_5_1} :
\begin{equation}
\label{Lap_5_2}
\tilde{\kappa}(p) =  \sum_{m,n=0}^{+\infty}\frac{(-\lambda_1)^m  (\gamma_1)_m (-\lambda_2)^n (\gamma_2)_n }{m! n!}p^{-\alpha_1 m - \alpha_2 n-\beta},
\end{equation}
\begin{equation}
\label{Lap_5_3}
\tilde{k}(p) = \sum_{m,n=0}^{+\infty}\frac{(-\lambda_1)^m  (-\gamma_1)_m (-\lambda_2)^n (-\gamma_2)_n }{m! n!}p^{-\alpha_1 m-\alpha_2 n+\beta-1},
\end{equation}
that are valid under the conditions $\alpha_1>0,\ \alpha_2>0,\ 0<\beta <1,\ \Re(p) > 0,\ |\lambda_1 p^{-\alpha_1}|<1,\ |\lambda_2 p^{-\alpha_2}|<1$.

To transform the functions  $\tilde{\kappa}$ and $\tilde{k}$ into the time domain, we again apply the formula \eqref{power_L} and  arrive at the representations
$$
\kappa(x)  = \sum_{m,n=0}^{+\infty}\frac{(-\lambda_1)^m  (\gamma_1)_m (-\lambda_2)^n (\gamma_2)_n }{m! n!}h_{\alpha_1 m+\alpha_2 n+\beta}(x) = 
$$
\begin{equation}
\label{time_5_4}
x^{\beta-1}\, \xi_2(\gamma_1;\gamma_2;(\alpha_1,\alpha_2;\beta);-\lambda_1x^{\alpha_1},-\lambda_2x^{\alpha_2}),
\end{equation}
$$
k(x)  =  \sum_{m,n=0}^{+\infty}\frac{(-\lambda_1)^m  (-\gamma_1)_m (-\lambda_2)^n (-\gamma_2)_n }{m! n!}h_{\alpha_1 m +\alpha_2 n +1-\beta}(x) = 
$$
\begin{equation}
\label{time_5_5}
x^{-\beta}\, \xi_2(-\gamma_1;-\gamma_2;(\alpha_1,\alpha_2;1-\beta);-\lambda_1x^{\alpha_1},-\lambda_2x^{\alpha_2}),
\end{equation}
where $\alpha_1>0,\ \alpha_2>0,\  0<\beta <1$ and $\xi_2$ is a new special function of the hypergeometric type in two variables defined by the convergent series
\begin{equation}
\label{X_3}
\xi_2(\gamma_1;\gamma_2;(\alpha_1,\alpha_2;\beta);y,z):= \sum_{m,n=0}^{+\infty} \frac{(\gamma_1)_m (\gamma_2)_n}{\Gamma(\alpha_1 m+\alpha_2 n+\beta)}\, \frac{y^mz^n}{m! n!},
\end{equation}
$$
\alpha_1>0,\ \alpha_2>0,\  \gamma_1,\gamma_2,\beta, y,z \in \Com. 
$$
\end{example}

Similar to the case considered in Example \ref{ex_4_4}, we denoted the function defined by the equation \eqref{X_3} by $\xi_2$ because for $\alpha_1 =\alpha_2 =1$  it reduces to the known  Horn function $\Xi_2$:
\begin{equation}
\label{Horn_x2}
 \xi_2(\gamma_1;\gamma_2;(1,1;\beta);y,z)= \sum_{m,n=0}^{+\infty} \frac{(\gamma_1)_m (\gamma_2)_n}{\Gamma(m+n+\beta)}\, \frac{y^mz^n}{m! n!}= \frac{1}{\Gamma(\beta)} 
\Xi_2(\gamma_1;\gamma_2;\beta;y,z),
\end{equation}
where the Horn function  $\Xi_2$ is  defined by the double confluent series of the hypergeometric type (see the formula (26) in section 5.7.1 in \cite{Bat}):
\begin{equation}
\label{Horn_x}
\Xi_2(\gamma_1;\gamma_2;\beta;y,z) := \sum_{m,n=0}^{+\infty} \frac{(\gamma_1)_m (\gamma_2)_n}{(\beta)_{m+n}}\, \frac{y^mz^n}{m! n!}.
\end{equation}

Thus, for $\alpha_1 = \alpha_2 = 1$, the pair \eqref{time_5_4} and \eqref{time_5_5} of the symmetrical Sonin kernels  is expressed in terms of the Horn function $\Xi_2$:
\begin{equation}
\label{time_5_6}
\kappa(x)  =  h_\beta(x)\, \Xi_2(\gamma_1;\gamma_2;\beta;-\lambda_1x ,-\lambda_2x ), \ \beta >0,
\end{equation}
\begin{equation}
\label{time_5_7}
k(x)  =  h_{1-\beta}(x)\, \Xi_2(-\gamma_1;-\gamma_2;1-\beta;-\lambda_1x ,-\lambda_2x), \ \beta <1.
\end{equation}

To the best knowledge of the author, both the pair of functions \eqref{time_5_4} and  \eqref{time_5_5} and  its particular case in form \eqref{time_5_6} and \eqref{time_5_7} are  new symmetrical Sonin kernels not yet reported in the literature. 

It is easy to see that all Sonin kernels that we considered in this paper, both the known and the new ones, belong to the space of functions $C_{-1}(0,+\infty)\, = \, \{\phi:\ \phi(x) = x^{p}\phi_1(x),\ x>0,\ p>-1,\ \phi_1\in C[0,+\infty)\}$. The GFIs and the GFDs with the Sonin kernels from this space were introduced and investigated in \cite{Luc_20}-\cite{Luc21c} and other related publications. 

In particular, the general theory of the GFIs and the GFDs developed in \cite{Luc_20}-\cite{Luc21c} ensures that on the suitable spaces of functions (see, e.g., \cite{Luc21a} for details), the GFD with the kernel $k$ as in \eqref{time_5_7}
$$
(\D_{(k)}\, f)(x) = \frac{d}{dx}\, \int_0^x h_{1-\beta}(x-\xi)\, \Xi_2(-\gamma_1;-\gamma_2;1-\beta;-\lambda_1(x-\xi) ,-\lambda_2(x-\xi))\, f(\xi)\, d\xi
$$
as well as the regularized GFD
$$
(\,_*\D_{(k)}\, f)(x) =  \int_0^x h_{1-\beta}(x-\xi)\, \Xi_2(-\gamma_1;-\gamma_2;1-\beta;-\lambda_1(x-\xi) ,-\lambda_2(x-\xi))\, f^\prime(\xi)\, d\xi
$$
are the left-inverse operators to the GFI with the kernel $\kappa$ as in \eqref{time_5_6}
$$
(\I_{(\kappa)}\, f)(x) = \int_0^x h_\beta(x-\xi)\, \Xi_2(\gamma_1;\gamma_2;\beta;-\lambda_1(x-\xi) ,-\lambda_2(x-\xi) )\, f(\xi)\, d\xi.
$$

Further properties of the new symmetrical Sonin kernels derived in this paper as well as  the GFIs and GFDs with these kernels will be considered elsewhere.

\end{document}